\documentclass[11pt]{amsart} 
\usepackage{amsthm}
\usepackage{amssymb, amsmath}
\usepackage{url,enumitem}

\theoremstyle{definition}

\numberwithin{equation}{section} 

\theoremstyle{remark}

\title[Rebuttal of Kowalenko's paper]{Rebuttal of Kowalenko's paper as concerns\\the irrationality of Euler's constant}
\author{Mark W. Coffey and Jonathan Sondow}
\address{Department of Physics, Colorado School of Mines, Golden, CO  80401, USA}
\email{mcoffey@mines.edu}
\address{209 West 97th Street \#6F, New York, NY  10025, USA}
\email{jsondow@alumni.princeton.edu}
\date{}         

\begin{document}


\newcommand{\CC}{\mathbb C} 


\newcommand{\beas}{\begin{eqnarray*}} 
\newcommand{\eeas}{\end{eqnarray*}} 

\newcommand{\bm}[1]{{\mbox{\boldmath $#1$}}} 

\newcommand{\bc}[2]{\genfrac{(}{)}{0pt}{}{#1}{#2}}

\newcommand{\tworow}[2]{\genfrac{}{}{0pt}{}{#1}{#2}}

\begin{abstract}
We rebut Kowalenko's claims in 2010 that he proved the irrationality of Euler's constant $\gamma$, and that his rational series for $\gamma$ is new.
\end{abstract}


\maketitle


\section{Introduction}

The irrationality of Euler's constant
$$\gamma=\lim_{n \to \infty}\left(\sum_{k=1}^n \frac{1}{k}-\log n\right)=0.577215664901532860606512090082402431042\dotso$$
has long been conjectured.  However, it remains an open problem.

In 2010 Kowalenko claimed that simple arguments suffice to settle this matter \cite{K}.  As
he offered no general framework or new mathematical principle, we believe that the following illustrations
are sufficient to describe the flaws in his very limited approach.

\section{A faulty irrationality argument}

Kowalenko derives the following formula for Euler's constant in equation (65) of \cite[p.~428]{K}:
\begin{equation*}
\gamma= \sum_{k=1}^{\infty}\frac{A_0}{k(k + 1)} - \sum_{k=2}^{\infty}\frac{A_1}{k^2} + \sum_{k=3}^{\infty}\frac{A_2}{k(k - 1)} - \sum_{k=4}^{\infty}\frac{A_3}{k(k - 2)}+\dotsb.
\end{equation*}
Here $A_0,A_1,A_2,A_3,\dotso$ are certain rational numbers. He writes:

\begin{quote}With the exception of the second series on the rhs, all the series in (65) can be easily evaluated
by decomposing them into partial fractions. On the other hand, the second series on the
rhs is virtually equal to $\zeta(2)$.
\end{quote}

He then transforms formula (65) into the following series in equation (69):
$$
\gamma= \frac32-\frac12 \left(\frac{\pi^2}{6} + \frac{1}{12} + \frac{5}{144} + \frac{247}{12960} +\frac{ 77}{6400} +\frac{25027}{3024000}+\dotsb\right).
$$
Kowalenko states:

\begin{quote}For $\gamma$ to be rational the term involving $\pi^2/6$, which arises solely from the summation
over $1/k^2$ or $\zeta(2)$ in (65), has to be cancelled by the remaining sum. This means that we need
to examine the methods for converting an irrational number into a rational number by the
process of addition. There are only two possible methods for achieving this, which are best
understood if we regard an irrational number as an infinite random distribution of decimal
digits. $\dotso$ Therefore, for $\gamma$ to be rational, we need to convert a random distribution
into a non-random one.
\end{quote}

But, for example, the distribution of digits in Liouville's irrational number
$$\sum_{n=1}^\infty \frac{1}{10^{n!}} =0.1100010000000000000000010 \ldots$$
is {\em not} random, as the sum formula shows. Thus Kowalenko's understanding of irrationality is lacking.

He continues:

\begin{quote}The first method by which an irrational number can be converted to a rational number is
to add another number, which at some stage possesses the opposite random distribution to
the original irrational number. This represents the situation whereby the second number can
be expressed as $C-\pi^2/6$, with $C$ a rational number. Such a situation, however, cannot occur
with (69). First, we note that if we are to obtain $C-\pi^2/6$ from the remaining terms after
the $\pi^2/6$ term in the parenthesis of (69), then these terms would have to yield a summation
involving $1/k^2$. This is simply not possible as all the $1/k^2$ terms have already been removed
as mentioned above.
\end{quote}

Here Kowalenko apparently assumes that Euler's series $\sum_{k=1}^\infty \frac{1}{k^2}=\frac{\pi^2}{6}$
is the only way to represent $\pi^2/6$ as the sum of a series of rational numbers. Of course, that is not true; there are infinitely many such representations of $\pi^2/6$, and hence of $C-\pi^2/6$, for any rational number $C>\pi^2/6$. For instance, taking $C=2$, an alternate series for $2-\pi^2/6$ is
$$2-\frac{\pi^2}{6} = \sum_{n=3}^{\infty}\sum_{k=2}^{\infty}\frac{1}{k^n}=\frac{1}{8}+\frac{1}{16}+\frac{1}{27}+\frac{1}{32}+\frac{1}{64}+\frac{1}{81}+\frac{1}{125}+\dotsb,$$
from Goldbach's theorem $\sum_{n=2}^{\infty}(\zeta(n)-1)=1$ (see, e.g., \cite[p.~142, equation (13)]{sc}).

Kowalenko concludes the paragraph:

\begin{quote}Furthermore, the summation would have to be negative. Yet all the
terms in the parenthesis in (69) are positive definite. Therefore, it is simply impossible for
all the terms in (69) to yield $C-\pi^2/6$.
\end{quote}

Here he claims that the sum of a series of positive rational numbers cannot be equal to $C-\pi^2/6$. But, for example, decimal expansion does give such a series:
$$C-\frac{\pi^2}{6} = n + 0.d_1d_2d_3\dotso = n+ \sum_{k=1}^{\infty}\frac{d_k}{10^k}.$$

In view of his misconceptions, Kowalenko has not proven that the ``first method by which an irrational number can be converted to a rational number'' does not lead to the rationality of Euler's constant. Its irrationality therefore remains an open problem.

\section{A known representation}

Finally, we point out that Kowalenko's claim to have found ``a new representation for Euler's constant'' \cite[p. 143]{K}
is also incorrect. Namely, what he calls ``Hurst's formula,'' which is $\gamma = \sum_{k=1}^\infty \frac{\lvert A_k\rvert}{k}$ in equation~(60) of \cite{K}, is known. According to Gourdon and Sebah \cite[p.~6]{GS}, the formula was discovered in 1924 by Kluyver \cite{Kl} (see the translation \cite{KlEng}), who wrote it as
$$
\gamma = \sum_{k=1}^\infty \frac{a_k}{k} = \frac{1}{2} + \frac{1}{24} + \frac{1}{72} + \frac{19}{2880} + \frac{3}{800} + \frac{868}{362880} + \frac{275}{169344} + \dotsb.
$$
By comparing the values of Kluyver's numbers $a_k$ (see the end of Section~2.3 in \cite[p.~6]{GS}) with those of Kowalenko's numbers $A_k$ (see Table~1 in \cite[p.~418]{K}) for $1\le k\le7$, one readily observes that $a_k = \lvert A_k\rvert$.

Indeed, letting
$$(z)_n=\frac{\Gamma(z+n)}{\Gamma(z)} = z(z+1)\dotsb(z+n-1)$$
denote
the Pochhammer symbol, with $\Gamma$ the Gamma function, we have the following relations.  From (18) of \cite{K},
$$A_k=\frac{(-1)^k}{k!} \int_0^1 (-x)_k dx,$$
while from \cite[p.~150]{Kl} and \cite[p. 143]{KlEng},
$$a_k=\frac{1}{k!}\int_0^1 x(1-x)_{k-1}dx=-\frac{1}{k!}\int_0^1 (-x)_k dx,$$
and we therefore conclude that $a_k=(-1)^{k+1}A_k$.

\section{Appendix: An anonymous referee's report}

The authors submitted a version of the first three sections of this paper to \emph{Acta Applicandae Mathematicae}, along with a list of six experts in irrationality theory as potential referees. The editors accepted the paper, writing in part, ``We have received the following clear and concise confirmation by one of the experts'':

\begin{quote}I completely agree with arguments of J.~Sondow and M.~Coffey. The article of V.~Kowalenko is baseless.\end{quote}


\end{document}